\documentclass[12pt]{amsart}

\newcommand{\bE}{\mathbb{E}}

\newcommand{\bP}{\mathbb{P}}
\newcommand{\bQ}{\mathbb{Q}}
\newcommand{\bR}{\mathbb{R}}

\newcommand{\bZ}{\mathbb{Z}}

\newcommand{\diag}{\mathrm{diag}}
\newcommand{\diam}{\mathrm{diam}}

\theoremstyle{plain} \newtheorem{Theo}{Theorem}[section]
\theoremstyle{plain} \newtheorem{Lemma}[Theo]{Lemma}
\theoremstyle{plain} 
\theoremstyle{plain} \newtheorem{Prop}[Theo]{Proposition}
\theoremstyle{remark} \newtheorem{Rem}[Theo]{Remark}
\theoremstyle{definition} 
\theoremstyle{definition} 
\theoremstyle{definition} \newtheorem{Ass}[Theo]{Assumption}
\theoremstyle{remark} 

\begin{document}

\title[Zeros of $p$-adic polynomials]{The expected number of zeros of a random system of $p$-adic polynomials}

\author{Steven N.\ Evans}
\address{Steven N.\ Evans \\
  Department of Statistics \#3860 \\
  University of California at Berkeley \\
  367 Evans Hall \\
  Berkeley, CA 94720-3860 \\
  U.S.A.}
\email{evans@stat.Berkeley.EDU}
\urladdr{http://www.stat.berkeley.edu/users/evans/}
\thanks{Supported in part by NSF grant DMS-0405778.  Part of the
research was conducted during a visit to the American Institute of Mathematics for a Workshop on Random Analytic Functions.}


\begin{abstract}
We study the simultaneous zeros of a random family of $d$ polynomials
in $d$ variables over the $p$-adic numbers.  For a family of natural
models, we obtain an explicit constant
for the expected number of zeros that lie
in the $d$-fold Cartesian product
of the $p$-adic integers. Considering models
in which the maximum degree that each variable
variable appears is $N$, this
expected value is
\[
p^{d \lfloor \log_p N \rfloor} \left(1 + p^{-1} + p^{-2} + \cdots + p^{-d}\right)^{-1}
\]
for the simplest such model.
\end{abstract}

\maketitle

\section{Introduction}

Various questions regarding the 
distribution of the number of real roots of
a random polynomial were considered in
 \cite{LitOff38, LitOff39, MR0009656}
and were taken up 
in \cite{MR0007812, MR0009655, MR0030713}, where
the main result is that the expected number of roots of a degree $n$ polynomial with independent standard Gaussian coefficients is
asymptotically equivalent to $\frac{2}{\pi} \log n$ for large $n$.
There has since been a huge amount of work on various aspects
of the distribution of the roots  of random polynomials 
and systems of random polynomials for
a wide range of models
with coefficients that are possibly dependent and have distributions
other than Gaussian.  It is impossible to survey this work adequately,
but some of the more commonly cited early papers are \cite{MR0234512, MR0234513, MR0286157, MR0288824}.  Reviews of the literature
can be found in \cite{MR856019, MR1290398, MR1376652, MR1679392}, and
some recent papers that indicate the level of sophistication
that has been achieved in terms of results and methodology
are \cite{MR1308023, MR1390040, MR1881723, MR1886266, 
MR1915821, MR1935565, MR2015330, MR2074895, MR2182444}.

In this paper we study the roots of random polynomials
over a field other than the real or complex numbers, the
field of $p$-adic numbers for some prime $p$.  
Like the reals, the $p$-adics arise
as a completion of the rationals with respect to certain metric -- see below.
They are the prototypical local fields (that is, non-discrete, locally
compact topological fields) and any local field with characteristic
zero is 
a finite algebraic extension of the $p$-adic numbers (the local
fields with non-zero characteristic are finite algebraic extensions
of the $p$-series field of Laurent series over the finite field
with $p$ elements).

In order to describe our results we need to give a little background.  For a fuller
treatment,
we refer the reader to 
\cite{MR791759} for an excellent introduction
to local fields and analysis on them.

We begin by defining the $p$-adic numbers.  Fix a positive prime $p$.
We can write any non-zero rational number $r \in \bQ \backslash \{ 0 
\}$ uniquely as $r=p^s (a/b)$
where $a$ and $b$ are not divisible by $p$. Set $| r | = p^{-s}$.
If we set $| 0 | =0$, then the map $| \cdot |$ has
the properties:
\begin{equation}
\label{valuation}
\begin{split}
| x | & = 0 \Leftrightarrow x=0, \\
| xy | & = | x | | y |, \\
| x+y | & \le | x | \vee | y |.
\end{split}
\end{equation}
The map $(x,y) \mapsto | x-y |$ defines a metric on $\bQ$, and we
denote the completion of $\bQ$ in this metric by $\bQ_p$.
The field operations on $\bQ$ extend continuously  to make $\bQ_p$
a  topological field called the {\em $p$-adic numbers}.
The map $| \cdot |$ also extends continuously and the extension
has properties (\ref{valuation}).  The closed unit ball around $0$,
$\bZ_p := \{x \in \bQ_p : |x| \le 1 \}$,
is the closure in $\bQ_p$ of the integers $\bZ$,
and is thus
a ring (this is also apparent from (\ref{valuation})),
called the  {\em $p$-adic integers}.
As
$\bZ_p = \{x \in \bQ_p : |x| < p\}$, the set $\bZ_p$ is also open.
Any other ball around $0$ is of the form
$\{x \in \bQ_p : |x| \le p^{-k}\} = p^k \bZ_p$
for some integer $k$.  Such a ball is the closure of the rational numbers
divisible by $p^k$, and is thus a $\bZ_p$-sub-module (this is again 
also apparent from (\ref{valuation})).  In particular, such a ball is an additive
subgroup of $\bQ_p$.  Arbitrary
balls are translates (=  cosets) of these closed and open subgroups.
In particular, the topology of $\bQ_p$ has a base of closed and open sets,
and hence $\bQ_p$ is totally disconnected.
Further, each of these balls is compact, and hence $\bQ_p$
is also locally compact.

There is a unique Borel measure $\lambda$
on $\bQ_p$ for which
\begin{equation*}
\begin{split}
\lambda (x+A) & = \lambda (A) , \quad x \in \bQ_p, \\
\lambda (xA) & = | x | \lambda (A), \quad x \in \bQ_p, \\
\lambda (\bZ_p) & =1. \\
\end{split}
\end{equation*}
The measure $\lambda$ is just suitably normalized Haar measure on the 
additive group of $\bQ_p$.
The restriction of $\lambda$ to $\bZ_p$ is 
the weak limit as $n \rightarrow \infty$ of the sequence of 
probability measures that at the $n$-th stage assigns mass $p^{-n}$ 
to each of the points  $\{0,1,\ldots,p^n-1\}$.

There is a substantial literature on
probability on the
$p$-adics and other local fields.
Two notable early papers are
\cite{MR772351, MR1111789}.
We have shown in a sequence papers \cite{MR990478, MR1118442,
MR1245397, MR1414928, MR1832433, MR1873668, MR1934156}
that the natural analogues on $\bQ_p$
of the centered Gaussian measures on $\bR$
are the normalized restrictions of $\lambda$ to the compact
$\bZ_p$-sub-modules $p^k \bZ_p$ and the point mass at $0$.  
More generally,
the natural counterparts of centered Gaussian measures for
$\bQ_p^d$ are normalized Haar measures on compact $\bZ_p$-sub-modules.
We call such probability measures {\em $\bQ_p$-Gaussian}
and say that a random variable distributed according to normalized
Haar measure on $\bZ_p^d$ is {\em standard} $\bQ_p$-Gaussian.
There are numerous papers on
Markov processes taking values in local fields, for example 
\cite{MR1163602, MR1290704, MR1705599, 
MR1893096, MR1778752, 
MR1822891, MR1900359, MR2034543, MR2156128}.  There are
also extensive surveys of the literature in the books
\cite{MR1746953, MR1848777, MR2105195}. 

If we equip the space of continuous
functions $C(\bZ_p^d, \bQ_p)$ with the map 
$f \mapsto \|f\| := \sup\{|f(t)| : t \in \bZ_p^d\}$, then
$\|\cdot\|$ is a $p$-adic norm in the sense that
\begin{equation*}
\begin{split}
\| f \| & = 0 \Leftrightarrow f=0, \\
\| a f \| & = | a | \| f \|, \quad a \in \bQ_p, \; f \in  C(\bZ_p^d, \bQ_p), \\
\| f+g \| & \le \| f \| \vee \| g \|.
\end{split}
\end{equation*}
Moreover, $C(\bZ_p^d, \bQ_p)$ is a $p$-adic Banach space in the sense
that it is complete with respect to the metric $(f,g) \mapsto \|f-g\|$.

There is a natural notion of orthogonality on the space 
$C(\bZ_p^d, \bQ_p)$. A collection $\{f_0, f_1, \ldots\}$
is {\em orthogonal} if $\|\sum_{k=0}^n a_k f_k\| = \bigvee_{k=0}^n |a_k| \|f_k\|$
for any $n$ and any $a_k \in \bQ_p$.  At first glance, this looks
completely unlike the notion of orthogonality one is familiar with
in real and complex Hilbert spaces, but it can be seen from 
\cite{MR791759} that there are actually close parallels.  
It is apparent from \cite{MR791759} that the sequence
of functions $\{t \mapsto \binom{t}{k}\}_{k=0}^\infty$, where
$\binom{t}{k} := \frac{t(t-1) \cdots (t-k+1)}{k!}$ (the
{\em Mahler basis}) is a very natural {\em orthonormal} basis
for $C(\bZ_p, \bQ_p)$ (that is, it is orthogonal and
each element has unit norm).  It is not hard to see that
the functions
\[
(t_1, t_2, \ldots, t_d) \mapsto 
\binom{t_1}{k_1} \binom{t_2}{k_2} \cdots \binom{t_d}{k_d},
\quad 0 \le k_1, k_2, \ldots, k_d < \infty,
\]
are an orthonormal basis for $C(\bZ_p^d, \bQ_p)$.

Putting all of these ingredients together, we see that
a natural model for a random system of $d$ independent
identically distributed $\bQ_p$-valued polynomials in $d$ variables lying
in $\bZ_p$ is the system
\[
F_i(t_1, t_2, \ldots, t_d) := \sum_k a_{k} Z_{i,k} \binom{t_1}{k_1} 
\binom{t_2}{k_2} \cdots \binom{t_d}{k_d}, \quad 1 \le i \le d,
\]
where the sum is over multi-indices $k=(k_1, k_2, \ldots, k_d)$,
for each $i$ the constants $a_{k} \in \bQ_p$ are zero for all but finitely
many $k$, and the $\bQ_p$-valued random variables are independent
and standard $\bQ_p$-Gaussian distributed.  

\begin{Ass}
\label{A:assumptions}
Assume that $a_{0} \ne 0$
and $a_{e_j} \ne 0$, $1 \le j \le d$,  where
$e_1 := (1,0,0, \ldots, 0)$, $e_2 := (0,1,0,\ldots, 0)$, and so on.
By re-scaling, we can assume without loss of generality that
$a_{0} = 1$.
We will also suppose that
$|a_{k}| \ge |a_{\ell}|$ when $k \le \ell$ in the usual
partial order on multi-indices (that is,
if $k=(k_1, k_2, \ldots, k_d)$ and
$\ell = (\ell_1, \ell_2, \ldots, \ell_d)$,
then $k_j \le \ell_j$ for $1 \le i \le d$).
It follows from the orthonormality of the products of Mahler basis elements
that each $(F_1, F_2, \ldots, F_d)$ maps $\bZ_p^d$
into $\bZ_p^d$.
\end{Ass}

\begin{Theo}
\label{T:main}
Suppose that Assumption \ref{A:assumptions}
holds.
For $(x_1, x_2, \ldots, x_d) \in \bZ_p^d$, the expected
number of points in the set 
\[
\{(t_1, t_2, \ldots, t_d) \in \bZ_p^d: F_i(t_1, t_2, \ldots, t_d) = x_i,
\; 1 \le i \le d\}
\]
is
\[
\left[\prod_{j=1}^d \bigvee_{h=1}^\infty \left|\frac{a_{h e_j}}{h}\right|\right]
\left(1 + p^{-1} + p^{-2} + \cdots + p^{-d}\right)^{-1}.
\]
\end{Theo}

Following some preliminaries
in Section~\ref{S:prelim}, 
we give the proof in Section~\ref{S:main}.  
However, we provide a
heuristic argument now as motivation for
the development we need to do in  Section~\ref{S:prelim}.  
Because we are arguing heuristically, we
do not justify various
interchanges of limits, sums and expectations.

Suppose first that $d=1$ and $x=0$.
Write $F = \sum_k a_k Z_k \binom{\cdot}{k_1}$ for $F_1$.  Let $B_{n,0},
B_{n,1}, \ldots, B_{n,p^n-1}$ be a list
of the balls of radius $p^{-n}$ in
$\bZ_p$, numbered so that $0 \in B_{n,0}$.
Let $I_{i,j}^{m,n}$ be the indicator of the event
that the graph of $F$ intersects 
$B_{m,i} \times B_{n,j}$. 
The number of zeros of $F$, 
$|\{t \in \bZ : F(t) = 0\}|$, is 
\[
\lim_{m \rightarrow \infty}
\lim_{n \rightarrow \infty}
\sum_i I_{i,0}^{m,n}.
\]
Because $Z_0$ is distributed according to
Haar measure on $\bZ_p$, the distribution of $z+F$
is the same for all $z \in \bZ_p$ and so the expectation
in question is also the expectation of
\[
\lim_{m \rightarrow \infty}
\lim_{n \rightarrow \infty}
p^{-n} \sum_i \sum_j I_{i,j}^{m,n}.
\]
As we observe in Section~\ref{S:main},
$F$ is a stationary process on $\bZ_p$
(this is not at all obvious and will hold
if and only if $|a_0| \ge |a_1| \ge \ldots$, hence
our assumption to this effect).
Consequently the expectation
in question is also the expectation of
\[
\lim_{m \rightarrow \infty}
\lim_{n \rightarrow \infty}
p^m p^{-n} \sum_j I_{0,j}^{m,n}.
\]

As in the real case,
polynomials look approximately like affine
functions on small scales, so for large $m$
the restriction of $F$ to the ball $B_{m,0}$
of radius $p^{-m}$ around $0$ is equivalent up to 
first order to a random affine function 
$t \mapsto Z_0 + b W t$,
where $W$ is standard $\bQ_p$-Gaussian and $b \in \bQ_p$
is a non-zero constant.
This implies first of all that for large $m$ 
the restriction
is injective, so that
$\lim_{n \rightarrow \infty} p^{-n} \sum_i I_{0,j}^{m,n}$
is the Haar measure of the image of $B_{m,0}$ by $F$.
Moreover, the image has Haar measure approximately
that of the image by $t \mapsto Z_0 + b W t$,
which is exactly $|b| |W| p^{-m}$.
Thus the expectation in question is nothing other than
the expectation of $|b| |W|$.  It remains to note
that $|W|$ takes the value $p^{-r}$ with
probability $p^{-r} - p^{-(r+1)}$ for $r=0,1,2,\ldots$
to conclude that the expectation of $|b W|$
is $|b| \sum_r (1-p^{-1}) p^{-2r} = |b|(1 + p^{-1})$.

Essentially the same heuristic argument works for 
general $d$. Once again the problem is reduced
to considering the expected Haar measure of the image
of a small ball by a random  affine function.  
Computing the actual
value of the expectation is more complicated
however, as it involves evaluating the expected value
of the determinant of the linear part of the affine
function.

This paper appears to be the first
to consider roots of random polynomials over the $p$-adic field.
There has been some work on
random polynomials over finite fields, see \cite{MR1142739, MR1230136, MR1404907, MR2092625, MR2097336}.

\section{Preliminaries}
\label{S:prelim}

Write $\lambda_d$ for the $d$-fold product measure
$\lambda^{\otimes d}$. Thus $\lambda_d$
is Haar measure on the additive group of $\bQ_p^d$
normalized so that $\lambda_d(\bZ_p^d) = 1$.
The Euclidean analogue of the following result
is well-known.

\begin{Lemma}
\label{L:linearchange}
For a Borel set $A \subseteq \bQ_p^d$ and a $d \times d$
matrix $H$, the set $H(A)$ has Haar measure 
$\lambda_d (H(A)) = |\det(H)| \lambda_d (A)$.
\end{Lemma}

\begin{proof} If $H$ is singular, then the
range of $H$ is a lower dimensional subspace of 
$\bQ_p^d$ and the result is obvious.  

Suppose then that $H$ is invertible.  
Write $GL(d,\bZ_p)$ for the space of $d \times d$
matrices that have entries in $\bZ_p$
and are invertible with the inverse also having entries
in $\bZ_p$.  By Cramer's rule, a matrix $W$ is in
$GL(d,\bZ_p)$ if and only if it has entries
in $\bZ_p$ and $|\det(W)|=1$. Moreover,
$GL(d,\bZ_p)$ is the set of linear isometries of 
$\bQ_p^d$ equipped with the metric derived from
the norm $|(x_1,x_2, \ldots, x_d)| = \bigvee_{i=1}^d |x_i|$
(see Section 3 of \cite{MR1934156}).
From the representation
of $H$ in terms of its elementary divisors, we have
\[
H = U \diag(p^{k_1}, p^{k_2}, \ldots, p^{k_d}) V,
\]
for integers $k_1, \ldots, k_d$ 
and matrices $U,V \in GL(d,\bZ_p)$ (see Theorem 3.1 of
\cite{MR1934156}).   Because
$|\det(U)| = |\det(V)|=1$, it follows that
$|\det(H)| = p^{-(k_1+ \cdots + k_d)}$.

From the uniqueness of Haar measure,
$\lambda_d  \circ U$ and $\lambda_d  \circ V$ are both constant
multiples of $\lambda_d $.  Both $U$ and $V$ 
map the ball $\bZ_p^d$
bijectively onto itself. 
Thus $\lambda_d  \circ U = \lambda_d  \circ V = \lambda_d$.  

Again from the uniqueness of Haar measure,
$\lambda_d  \circ \diag(p^{k_1}, p^{k_2}, \ldots, p^{k_d})$
is a constant multiple of $\lambda_d$.  Now
\[
\begin{split}
& \lambda_d  \circ \diag(p^{k_1}, p^{k_2}, \ldots, p^{k_d})(\bZ_p^d)
 =
\lambda_d \left(\prod_{j=1}^d p^{k_j} \bZ_p \right) \\
& \quad =
\prod_{j=1}^d \lambda(p^{k_j} \bZ_p) 
 =
p^{-(k_1+ \cdots + k_d)} \\
& \quad =
|\det(H)|
 =
|\det(H)| \lambda_d(\bZ_p^d). \\
\end{split}
\]
\end{proof}

Write $gl(d, \bQ_p)$ for the space of $d \times d$ matrices
with entries in $\bQ_p$.
We say that a function $f$ from an open subset $X$ of $\bQ_p^d$
into $\bQ_p^d$ is {\em continuously differentiable}
if there exists a continuous function $R: X \times X \rightarrow gl(d, \bQ_p)$
such that $f(x)-f(y) = R(x,y)(x-y)$.  This definition is a natural generalization
of Definition 27.1 of \cite{MR791759} for the case $d=1$.  Set $J f(x) = R(x,x)$.

The next result is along the lines of the Euclidean
implicit function theorem.  It follows from 
Lemma~\ref{L:linearchange} and
arguments similar to those which establish the analogous
results for $d=1$ in Proposition 27.3, Lemma 27.4,
and Theorem 27.5 of
\cite{MR791759}.

\begin{Lemma}
\label{L:implicit}
Suppose for some open subset $X$ of $\bQ_p^d$
that $f:X \rightarrow \bQ_p^d$ is continuously
differentiable.
\begin{itemize}
\item[(i)]
If $J f(x_0)$ is 
invertible for some $x_0 \in X$, then, for all
sufficiently small balls $B$ containing $x_0$, 
the function
$f$ restricted to $B$ is a bijection onto its image,
$f(B) = J f(x_0) (B)$, and $|\det(J f(x))| = |\det(J f(x_0))|$
for $x \in B$.  In particular,
\[
\lambda_d (f(B)) = |\det(J f(x_0))| \, \lambda_d (B).
\]
\item[(ii)]
If $J f(x_0)$ is singular for some $x_0 \in X$, then,
for all sufficiently small balls $B$ containing $x_0$,
$\lambda_d (f(B)) = {\mathrm{o}}(\lambda_d (B))$.
\end{itemize}
\end{Lemma}

The following result is an analogue of a particular
instance of Federer's co-area formula.  The special
case of this result for $d=1$ and an injective
function is the substitution formula in Appendix A.7
of \cite{MR791759}.

\begin{Prop}
\label{P:coarea}
Suppose for some open subset $X$ of $\bQ_p^d$
that $f:X \rightarrow \bQ_p^d$ is continuously
differentiable.  Then, for any non-negative Borel function 
$g: \bQ_p^d \rightarrow \bQ_p^d$,
\[
\int_X g \circ f(x) \, |\det(J f(x))| \, \lambda_d (dx)
=
\int_{\bQ_p^d} g(y) \, \# f^{-1}(y) \, \lambda_d (dy).
\]
\end{Prop}

\begin{proof}
It suffices to consider the case when $g$ is the indicator
function of a ball $C$.  Write $\delta$ for the diameter of $C$.
Put
\[
S := \{x \in X: J f(x) \text{ is singular}\}
\]
and
\[
I := \{x \in X: J f(x) \text{ is invertible}\}.
\]

From Lemma~\ref{L:implicit}(ii), $\lambda_d (f(S))=0$,
so that 
\[
\lambda_d (\{y \in \bQ_p^d : f^{-1}(y) \cap S \ne \emptyset\}) = 0
\]
and 
\[
\begin{split}
\int_{\bQ_p^d} g(y) \, \# (f^{-1}(y) \cap S) \, \lambda_d (dy)
& =
0 \\
& =
\int_S g \circ f(x) \, |\det(J f(x))| \, \lambda_d (dx). \\
\end{split}
\]

From Lemma~\ref{L:implicit}(iii), we can cover the open set $I$
with a countable collection of balls $B_k$ such that
$f$ restricted to $B_k$ is a bijection onto its image,
$f(B) = J f(x_0) (B)$ for some
$x_0 \in B$, $|\det(J f(x))| = |\det(J f(x_0))|$
for all $x \in B_k$, $\lambda_d (f(B_k)) = |\det(J f(x_0))| \lambda_d (B_k)$, and $\diam f(B_k) \le \delta$, so that
$g$ is constant on $f(B_k)$.  Hence
\[
\begin{split}
& \int_{\bQ_p^d} g(y)  \, \# (f^{-1}(y) \cap B_k) \, \lambda_d (dy) \\
& \quad =
\int_{f(B_k)} g(y)  \, \lambda_d (dy) \\
& \quad =
\int_{B_k} g \circ f(x)\,  |\det(J f(x))| \, \lambda_d (dx) \\
\end{split}
\]
Summing over $k$ gives 
\[
\int_{\bQ_p^d} g(y)  \, \# (f^{-1}(y) \cap I) \, \lambda_d (dy)
=
\int_{I} g \circ f(x) \,  |\det(J f(x))| \, \lambda_d (dx)
\]
and the result follows.
\end{proof}

\section{Proof of Theorem~\ref{T:main}}
\label{S:main}

For $x \in \bZ_p^d$, write $N(x)$ for the
number of points in the set 
\[
\{(t_1, t_2, \ldots, t_d) \in \bZ_p^d: F_i(t_1, t_2, \ldots, t_d) = x_i,
\; 1 \le i \le d\}.
\]
Since $Z_{i,0} - (x_1,x_2,\ldots,x_d)$ has the same
distribution as $Z_{i,0}$, it follows that
$\bE[N(\cdot)]$ is constant.
Also,
by an extension of the argument for $d=1$ in Theorem 9.3
of \cite{MR990478} (see also Theorem 8.2 of 
\cite{MR1832433}), 
the stochastic processes $F_i$ are stationary.

Thus, by Proposition~\ref{P:coarea},
\[
\begin{split}
\bE[N(x)]
& =
\int_{\bZ_p^d} \bE[N(x)] \, \lambda_d(dx) \\
& =
\bE\left[\int_{\bZ_p^d} N(x) \, \lambda_d(dx) \right] \\
& =
\bE\left[\int_{\bZ_p^d} |\det(J F(t))|  \, \lambda_d(dt)\right] \\
& =
\int_{\bZ_p^d} \bE[|\det(J F(t))|]  \, \lambda_d(dt) \\
& =
\bE[|\det(J F(0))|]. \\
\end{split}
\]

Now
\[
(J F(0))_{ij} = \sum_h a_{h e_j} Z_{i, h e_j} \frac{(0-1)(0-2) \ldots (0-h+1)}{h!}
= b_j W_{ij},
\]
where the $W_{ij}$ are standard $\bQ_p$-Gaussian random variables and
$b_j \in \bQ$ is any constant with 
\[
|b_j| = \bigvee_h \left|\frac{a_{h e_j}}{h}\right|,
\]
and so 
\[
\det(J F(0))
=
\left[\prod_{j=1}^d b_j \right]  \det\left(W_{ij} \right)_{1 \le i,j \le d}.
\]

From Theorem 4.1 in \cite{MR1934156}, 
we find, putting
\[
\Pi_k := (1 - p^{-1})(1 - p^{-2}) \cdots (1 - p^{-k}),
\]
that 
\[
\begin{split}
\bE[|\det(J F(0))|]
& =
\left[\prod_{j=1}^d |b_j|\right]
\sum_{h=0}^\infty p^{-h} \bP\{|\det\left(W_{ij} \right)_{1 \le i,j \le d}| = p^{-h}\} \\
& =
\left[\prod_{j=1}^d |b_j|\right]
\sum_{h=0}^\infty
p^{-2h} \frac{\Pi_d \Pi_{d+h-1}}{\Pi_h \Pi_{d-1}}. \\
\end{split}
\]
The result then follows from a consequence of the
$q$-binomial theorem, see Corollary 10.2.2 of \cite{MR1688958}.

\begin{Rem}
(i) Suppose that $a_{(k_1, \ldots, k_d)} = 1$ if $k_i \le N$ for all $i$
and is zero otherwise.  Then $|b_j|$ is just $p^r$, where $r = \lfloor \log_p N \rfloor$
is the largest power of $p$ that divides some integer $\ell$ with $1 \le \ell \le N$.

\noindent
(ii) Results about level sets of Euclidean processes
are often obtained using the Kac-Rice formula.
As shown in \cite{MR2115043}, result like the
Kac-Rice formula are a consequence of Federer's
co-area formula (see also \cite{AdTa06} for
an extensive discussion of this topic).
It would be possible to derive
a $p$-adic analogue of the Kac-Rice
formula from Proposition~\ref{P:coarea} and use it to
prove Theorem~\ref{T:main}.  However, the homogeneity
in ``space'' of $(F_1, F_2, \ldots, F_d)$ makes this
unnecessary.  

\noindent
(iii) Because  $(F_1, F_2, \ldots, F_d)$ is stationary,
its level sets are all stationary point processes
on $\bZ_p^d$ with intensity the multiple of
$\lambda_d$ given in Theorem~\ref{T:main}.

\noindent
(iv) The requirement that the $F_i$ are identically distributed
could be weakened.  All we actually use is that the distribution of
$(J F(0))_{ij}$ does not depend on $i$.
\end{Rem}

\noindent
{\bf Acknowledgment:}  We thank two anonymous referees for suggestions that
improved the presentation of the paper.

\def\polhk#1{\setbox0=\hbox{#1}{\ooalign{\hidewidth
  \lower1.5ex\hbox{`}\hidewidth\crcr\unhbox0}}}
  \def\polhk#1{\setbox0=\hbox{#1}{\ooalign{\hidewidth
  \lower1.5ex\hbox{`}\hidewidth\crcr\unhbox0}}}
  \def\polhk#1{\setbox0=\hbox{#1}{\ooalign{\hidewidth
  \lower1.5ex\hbox{`}\hidewidth\crcr\unhbox0}}}
  \def\polhk#1{\setbox0=\hbox{#1}{\ooalign{\hidewidth
  \lower1.5ex\hbox{`}\hidewidth\crcr\unhbox0}}}
  \def\polhk#1{\setbox0=\hbox{#1}{\ooalign{\hidewidth
  \lower1.5ex\hbox{`}\hidewidth\crcr\unhbox0}}}
  \def\polhk#1{\setbox0=\hbox{#1}{\ooalign{\hidewidth
  \lower1.5ex\hbox{`}\hidewidth\crcr\unhbox0}}}
\providecommand{\bysame}{\leavevmode\hbox to3em{\hrulefill}\thinspace}
\providecommand{\MR}{\relax\ifhmode\unskip\space\fi MR }
\providecommand{\MRhref}[2]{%
  \href{http://www.ams.org/mathscinet-getitem?mr=#1}{#2}
}
\providecommand{\href}[2]{#2}

\end{document}